\documentclass[a4paper]{amsart}

\usepackage{amsmath}
\usepackage{amsfonts}
\usepackage{amsthm}
\usepackage{amsopn}
\usepackage{amssymb}
\usepackage{amscd}
\usepackage{ae}
\usepackage[all]{xy}
\DeclareMathOperator{\sign}{sign} \DeclareMathOperator{\ind}{ind}

\DeclareMathOperator{\res}{res} 
\DeclareMathOperator{\mo}{mod}

\DeclareMathOperator{\ann}{ann}
\DeclareMathOperator{\tor}{Tor}

\usepackage{amsthm}                 
\usepackage{amssymb}

\usepackage{mathrsfs}
\usepackage{dsfont}
\newcommand{\signa}{\textnormal{signature}}

\newcommand{\CC}{{\mathds{C}}}

\newcommand{\RR}{{\mathds{R}}}
\newcommand{\ZZ}{{\mathds{Z}}}

\newcommand{\OCn}{\mathscr{O}_{\CC^n,0}}

\newcommand{\ICn}{\ind_{\CC^n,0}}
\newcommand{\IVn}{\ind_{V,0}}

\newcommand{\RCn}{\res_{\CC^n,0}}

\newcommand{\RVn}{\res_{V,0}}

\newcommand{\dep}{\textnormal{depth}}
\newcommand{\ERn}{\mathscr{E}_{\RR^n,0}}
\newcommand{\IRn}{\ind_{\RR^n,0}}

\newcommand{\bil}{<,>_l}

\newcommand{\IVRn}{\ind_{V^{\RR},0}}

\newtheorem{thm}{Theorem}[section]
\newtheorem{satz}[thm]{Proposition}
\newtheorem{lemma}[thm]{Lemma}

\theoremstyle{definition}
\newtheorem{defn}[thm]{Definition}

\begin{document}

\title[Real and complex indices of vector fields on curves]
{Real and complex indices of vector fields on complete
intersection curves with isolated singularity}
\author{Oliver Klehn}
\email{klehn@math.uni-hannover.de}
\address{Institut für Mathematik, Universität Hannover,
        Postfach 6009, D-30060 Hannover, Germany }

\subjclass[2000]{Primary 32S65; Secondary 14B05, 13H10.}
\keywords{index, vector field, local residue, signature, socle, isolated singularity.}

\begin{abstract}
If $(V,0)$ is an isolated complete intersection singularity and $X$ a holomorphic vector field tangent
to $V$ one can define an index of $X$, the so called GSV index, which generalizes the Poincar\'{e}-Hopf index.
We prove that the GSV index coincides with the dimension of a certain explicitely constructed vector space, if
$X$ is deformable in a certain sense and $V$ is a curve. We also give a sufficient algebraic criterion for $X$ to be deformable
in this way. If one considers the real analytic case one can also define an index of $X$ which is called the
real GSV index. Under the condition that $X$ has the deformation property, we prove a signature formula for the
index generalizing the Eisenbud-Levine Theorem.
\end{abstract}

\maketitle

\section{Introduction}
\label{sec:introduction}
\subsection{Classical results}
Assume that the continuous map germ $g\colon (\RR^n,0)\to (\RR^n,0)$ defines an
isolated zero. Then the map
$g/||g||\colon S_{\delta}^{n-1}\to S^{n-1}$ of spheres around the origin has a degree, the so called
Poincar\'{e}-Hopf index $\IRn (g)$ of $g$. If $g$ is analytic one has algebraic interpretations of
this index, that we first want to describe. If $g\colon (\CC^n,0)\to (\CC^n,0)$ is holomorphic let
$Q_g$ be the algebra obtained by factoring $\OCn$ by the ideal generated by the components of $g$.
One has
\begin{thm}[\cite{agv, gh}]
\[
\ICn (g)=\dim_{\CC}Q_g.
\]
\end{thm}
Here we have made the identification $\CC^n\cong\RR^{2n}$ of course.
Now let $\ERn$ be the ring of real analytic function germs on $(\RR^n,0)$ and further
$g\colon (\RR^n,0)\to (\RR^n,0)$ be finite and real analytic, in the sense that
$Q_g$ is finite dimensional as $\RR$-vector space and where $Q_g$ is the algebra obtained by factoring
$\ERn$ with the ideal generated by the components of $g$ in this case.
If one denotes by $J_g$ the determinant of the Jacobian of $g$ one has the following famous
theorem:

\begin{thm}[Eisenbud-Levine]
\label{elt}
Let $l\colon Q_g\to \RR$ be a linear form with $l(J_g)>0$. Then
\[
\IRn (g)=\signa\bil .
\]
\end{thm}
Here $\bil$ is the induced bilinear form defined by $<h_1, h_2>_l:= l(h_1\cdot h_2)$.

\subsection{Generalization to complete intersections}
Now let $(V,0):=(\{ f_1=\dots =f_q=0\},0)\subset (\CC^n,0)$ be an isolated singularity of a complete intersection (ICIS) and
$X:=\sum_{i=1}^n X_i\frac{\partial}{\partial z_i}$
be the germ of a holomorphic vector field on $(\CC^n,0)$ tangent to $V$, say $Xf=Cf$ with an isolated
zero on $V$. In this situation one can also define an index $\IVn (X)$, called the (complex) GSV index,
see \cite{asv, bg, gsv}, and it is the Poincar\'{e}-Hopf index when $V$ is smooth. The definition of the index is as follows:

Choose a sufficiently small sphere $S_{\delta}$ around the origin in $\CC^n$ which intersects $V$ transversally
and consider the link $K=V\cap S_{\delta}$ of $V$. The vectors $X,\nabla f_1,\dots ,\nabla f_q$ are linearly independent for all points
of $K$ and we have a well defined map
\[
(X, \nabla f_1,\dots ,\nabla f_q)\colon K\to W_{q+1}(\CC^n)
\]
where $W_{q+1}(\CC^n)$ denotes the manifold of $(q+1)$-frames in $\CC^n$ and we consider the complex gradients of course. We have
\[
H_{2n-2q-1}(K)\cong\ZZ,\;\;H_{2n-2q-1}(W_{q+1}(\CC^n))\cong\ZZ
\]
and therefore the map has a degree. We let $K$ to be oriented as boundary of the complex manifold $V\setminus\{ 0\}$ here.
The index $\IVn (X)$ of $X$ is defined to be the degree of this map. (If $V$ is a curve $K$ can have more components, we will sum
over the degrees of the components then.)

We now want to formulate our main theorems. We need a definition first.

\begin{defn}
$X$ is called a good vector field (wrt. $V$), if there is a holomorphic deformation $X_t$ of $X$,
so that for all $t\in\CC^q$ sufficiently close to zero $X_t$ is tangent to the
$t$-fibre $V_t$ of $f$. $X_t$ is called a good deformation of $X$.
\end{defn}
We will prove a sufficient criterion for a vector field to be good, which states that $X$ is good whenever
all coefficients of the matrix $C$ are contained in the ideal generated by the maximal minors of the Jacobian of $f$
in $\OCn$. It follows from the definition of the index that it equals the sum of the indices
of a good deformation on a smooth fibre.

After a linear generic change of coordinates one can assume that
$(f_1,\dots ,f_q, X_1,\dots X_{n-q})$
is a regular $\OCn$-sequence, see \cite{lss},
and we always assume the coordinates to be chosen in this way in this article.
Let $\mathscr{B}_0:=\OCn /(f_1,\dots ,f_q,X_1,\dots ,X_{n-q})$. Due to the chosen coordinates $\mathscr{B}_0$
is finite dimensional as complex vector space.
We also set $\mathscr{C}_0:=\mathscr{B}_0/\ann_{\mathscr{B}_0}(DF)$, where
\[
DF:=\det \left(\frac{\partial (f_1,\dots ,f_q)}{\partial (z_{n-q+1},\dots ,z_n)}\right).
\]

We prove an index formula for vector fields in the case $q=n-1$:
\begin{thm}
\label{them1}
Let $X$ be a good vector field and $V$ a curve. Then
\[
\IVn (X)=\dim_{\CC}\mathscr{C}_0.
\]
\end{thm}

Now let
\[
(V^{\RR},0):=(\{ f^{\RR}_ 1=\dots =f^{\RR}_ q =0\} ,0)\subset (\RR^n,0)
\]
defined by real analytic function germs. If $f$ denotes the complexification of $f^{\RR}$ we assume
that $f$ defines an ICIS of dimension $n-q$.
Furthermore let the real analytic vector field $X^{\RR}$ be tangent to $(V_{\RR},0)$ with an algebraic
isolated zero on $(V^{\RR},0)$.
One defines the real GSV index of $X^{\RR}$ similarly to the complex index, see \cite{asv},
and denotes this index by $\IVRn (X^{\RR})$ if $n-q$ is odd and by $\IVRn^2 (X^{\RR})$ if $n-q$ is even.
Due to topological reasons one can only define an $(\mo 2)$-index if $n-q$ is even. The definition of $X^{\RR}$
to be good is as in the complex case using real analytic deformations.

We prove for the case $q=n-1$ the following formula generalizing the Eisenbud-Levine theorem:

\begin{thm}
\label{them2}
Let $V^{\RR}$ be a curve, $X^{\RR}$ a good vector field and $l\colon\mathscr{C}_0^{\RR}\to\RR$ a linear form with $l(c_1)>0$.
Then
\[
\ind_{V^{\RR},0}(X^{\RR})=\signa\bil.
\]
\end{thm}
Here $\mathscr{C}_0^{\RR}$ is defined as in the complex case using $\ERn$ instead of $\OCn$,
$c_1$ is the coefficient of $t$ in the formal power series expansion
of $\det (\mathds{1}+tDX^{\RR})/\det (\mathds{1} +tC)$, where $DX^{\RR}$ is the Jacobian of $X^{\RR}$, and
$\bil$ is the induced bilinear form defined as in the classical case. $C$ is defined by the tangency
condition $X^{\RR}f^{\RR}=Cf^{\RR}$. 

\section{Residues of holomorphic vector fields}
To prove our main theorems we need a few results on residues of holomorphic vector fields that we
want to collect in this section.

Let $g\colon (\CC^n,0)\to (\CC^n,0)$ be a holomorphic map germ with isolated zero. Then
the residue
$\RCn^g(h)$ of any $h\in\OCn$ wrt. $g$ is defined as
\[
\RCn^g(h):=\frac{1}{(2\pi i)^n}\int_{\Gamma}\frac{hdz_1\wedge\dots\wedge dz_n}{g_1\cdot\ldots\cdot g_n},
\]
where $\Gamma$ is the real $n$-cycle
$\Gamma :=\{ |g_i|=\epsilon_i ,i=1,\dots ,n\}$ for $\epsilon_i\in\RR_{>0}$ small enough with orientation given by
$d( \arg g_1)\wedge\dots\wedge d(\arg g_n)\geq 0$.
Sometimes we denote this residue also by
\[
\RCn\begin{bmatrix} h\\g_1\dots g_n\end{bmatrix}.
\]
If we denote by $J_g$ the Jacobian determinant of $g$ one has the following classical result:
\begin{thm}[\cite{agv, gh}]
\label{res1}
\begin{itemize}
\item[(i)] $\RCn^g\colon Q_g\to\CC$ defines a linear form.
\item[(ii)] The induced bilinear form $<,>_{\RCn^g}$ is non-degenerate.
\item[(iii)] $\ICn (g)=\dim_{\CC}Q_g=\RCn^g(J_g)$.
\end{itemize}
\end{thm}
If we consider linear forms $l\colon Q\to\mathds{F}$ on commutative $\mathds{F}$-algebras for an arbitrary field
in this article, the induced bilinear form $\bil$ on $Q$ is always the bilinear form defined by
$<h_1,h_2>_l:=l(h_1\cdot h_2)$. The second statement in the theorem is usually called  "Local (Grothendieck-)Duality" and it
states that $Q_g$ is a Gorenstein algebra. This means that the annihilator of the maximal ideal, the socle, of $Q_g$ is
one-dimensional and it is well known that it is generated by the class of $J_g$. One immediately concludes
that for any linear form $l\colon Q_g\to\CC$ with $l(J_g)\neq 0$ the induced pairing $\bil$ is
non-degenerate.

As in the second part of the introduction
let $(V,0):=(\{ f_1=\dots =f_q=0\},0)\subset (\CC^n,0)$ be an isolated singularity of a complete intersection (ICIS) and
$X:=\sum_{i=1}^n X_i\frac{\partial}{\partial z_i}$
be the germ of a holomorphic vector field on $(\CC^n,0)$ tangent to $V$, say $Xf=Cf$ with an isolated
zero on $V$. Further let
\[
\Sigma :=\{f_1=\dots =f_q=0, |X_1|=\epsilon_1,\dots ,|X_{n-q}|=\epsilon_{n-q}\}
\]
be a small real $(n-q)$-cycle with orientation determined by
\[
d(\arg X_1)\wedge\dots\wedge d(\arg X_{n-q})\geq 0.
\]
Then we define the relative residue of any $h\in\OCn$ wrt. $X$ to be
\[
\RVn^X(h):=\frac{1}{(2\pi i)^{n-q}}\int_{\Sigma}\frac{hdz_1\wedge\dots\wedge dz_{n-q}}
{X_1\cdot\ldots\cdot X_{n-q}}.
\]
The absolute residue of $h$ wrt. $X$ is defined as
\[
\RCn^X(h):=\RCn\begin{bmatrix} h\\ X_1\dots X_{n-q} f_1\dots f_q\end{bmatrix}.
\]
Now let $c_{n-q}$ be the coefficient of $t^{n-q}$ in the formal power series expansion
of $\det (\mathds{1}+tDX)/\det (\mathds{1} +tC)$, where $DX$ is the Jacobian of $X$. We have the
following theorem proven in \cite{lss}:
\begin{thm}[Lehmann, Soares, Suwa]
\[
\IVn (X)=\RVn^X(c_{n-q})
\]
\end{thm}
The author has proven in \cite{k} that one always has $\RVn^X(h)=\RCn^X(hDF)$, where we have set
\[
DF:=\det \left(\frac{\partial (f_1,\dots ,f_q)}{\partial (z_{n-q+1},\dots ,z_n)}\right).
\]
This means that we also have $\IVn (X)=\RCn^X(c_{n-q}DF)$, which is one of the main tools in the proof of
our main theorems.

\section{Vector fields tangent to smooth varieties}
To prove our main theorems we first prove them for the smooth situation, what is done in this section,
and use good deformations
to generalize to the singular case. We also look at the socle of $\mathscr{C}_0$.
We use the notations as introduced in the first and second section.

\subsection{The complex situation}
Let $1\leq i_1<\dots <i_q\leq n$ be fixed and $1\leq j_1<\dots <j_{n-q}\leq n$
the complement of $i_1,\dots ,i_q$ in $\{1,\dots ,n\}$. Furthermore we set
\[
DI:=\det \left(\frac{\partial (f_1,\dots ,f_q)}{\partial (z_{i_1},\dots ,z_{i_q})}\right)
\]
and let $\sigma_I$ be the permutation defined by
\[
\sigma_I :=\begin{pmatrix}
1 & \dots & n-q & n-q+1 & \dots & n\\ j_1 & \dots & j_{n-q} & i_1 & \dots & i_q
\end{pmatrix},
\]
where $I$ is the multiindex $I:=(i_1,\dots ,i_q)$.
We also set
\[
\mathscr{B}^I_0:=\frac{\OCn }{(f_1,\dots ,f_q,X_{j_1},\dots ,X_{j_{n-q}})}
\]
and if $DI(0)\neq 0$
\[
\gamma_I :=\frac{\sign\;\sigma_I}{DI}\det\left( \frac{\partial (X_{j_1},\dots ,X_{j_{n-q}},f_1,\dots ,f_q)}
{\partial (z_1,\dots ,z_n)}\right) .
\]
\begin{lemma}
\label{lem1}
Let $DI(0)\neq 0$. Then
\[
\IVn (X)=\dim_{\CC}\mathscr{B}^I_0.
\]
\end{lemma}

\begin{proof}
By the implicit mapping theorem it is not hard to show, that
$X|_V$ corresponds to the vector field
\[
Y:=X_{j_1}\circ\psi\frac{\partial}{\partial y_1}+\dots +
X_{j_{n-q}}\circ\psi\frac{\partial}{\partial y_{n-q}}
\]
on $(\CC^{n-q},0)$, where $\psi\colon (\CC^{n-q},0)\to (V,0)$ is a biholomorphic map as in the implicit mapping
theorem.
From
\[
\frac{\mathscr{O}_{\CC^{n-q},0}}{(X_{j_1}\circ\psi ,\dots ,X_{j_{n-q}}\circ\psi )}\cong
\mathscr{B}^I_0
\]
the claim follows.
\end{proof}

\begin{lemma}
\label{lem22}
Let $DI(0)\neq 0$. Then we have for any $h\in\OCn$
\[
\RCn^X(hDF)=\sign\sigma_I\RCn\begin{bmatrix}
hDI\\ X_{j_1}\dots X_{j_{n-q}}f_1\dots f_q\end{bmatrix}.
\]
\end{lemma}

\begin{proof}
By the transformation formula for residues, see \cite{gh}, we have to show that there is a matrix $A$ with
\[
(X_1,\dots , X_{n-q}, f_1,\dots , f_q)^t=
A(X_{j_1},\dots ,X_{j_{n-q}},f_1,\dots ,f_q)^t
\]
and $\det A=\sign\;\sigma_I DF/DI$. From $Xf=Cf$ we get
\begin{equation}
\label{rech1}
\begin{pmatrix}X_{i_1}\\ \vdots \\ X_{i_q}\end{pmatrix}=Cf-
\left( \frac{\partial (f_1,\dots ,f_q)}{\partial (z_{i_1},\dots ,z_{i_q})}\right)^{-1}
\frac{\partial (f_1,\dots ,f_q)}{\partial (z_{j_1},\dots ,z_{j_{n-q}})}
\begin{pmatrix}X_{j_1}\\ \vdots \\ X_{j_{n-q}}\end{pmatrix}.
\end{equation}
Now let $i_1,\dots ,i_k\in\{1,\dots ,n-q\}$ and $i_{k+1},\dots ,i_q\in\{n-q+1,\dots ,n\}$.
Then it follows that $j_1,\dots ,j_{n-q-k}\in\{1,\dots ,n-q\}$ and
$j_{n-q-k+1},\dots ,j_{n-q}\in\{n-q+1,\dots ,n\}$.
For $k=0$ the claim follows immediately. With equation \ref{rech1} we obtain a matrix $B$ with
\[
(X_{j_1},\dots ,X_{j_{n-q+k}},X_{i_1},\dots ,X_{i_k},f_1,\dots ,f_q)^t=B(X_{j_1},\dots ,X_{j_{n-q}},f_1,\dots ,f_q)^t.
\]
If $\sigma'\in S_{n-q}$ is the permutation with
\[
\sigma'=\begin{pmatrix}1 & \dots & n-q+k & n-q+k+1 & \dots & n-q\\
j_1 & \dots & j_{n-q-k} & i_1 & \dots & i_k\end{pmatrix},
\]
we have $\det A=\sign\,\sigma'\det B$.
Then $\det B$ is the determinant of the upper right $(k\times k)$-block of the matrix
\[
-
\left( \frac{\partial (f_1,\dots ,f_q)}{\partial (z_{i_1},\dots ,z_{i_q})}\right)^{-1}
\frac{\partial (f_1,\dots ,f_q)}{\partial (z_{j_1},\dots ,z_{j_{n-q}})}.
\]
If $d_{l,m}$ is the determinant of the matrix obtained by replacing the $l$-th column
of
\[
\frac{\partial (f_1,\dots ,f_q)}{\partial (z_{i_1},\dots ,z_{i_q})}
\]
by the $m$-th column of
\[
\frac{\partial (f_1,\dots ,f_q)}{\partial (z_{j_1},\dots ,z_{j_{n-q}})}
\]
and if
\[
D:=(d_{l,m})_{l=1,\dots ,k}^{m=n-q-k+1,\dots ,n-q},
\]
we have to show
\[
\det D=(DI)^{k-1}
\det\left( \frac{\partial (f_1,\dots ,f_q)}{\partial (z_{j_{n-q-k+1}},\dots ,z_{j_{n-q}},
z_{i_{k+1}},\dots ,z_{i_q})}\right).
\]
This can be done by induction over $k$ where $k=1$ is obviously.
The conclusion is straightforward and we do not want to write it down here.
\end{proof}

\begin{lemma}
\label{lem3}
If $V$ is smooth, then
\[
\IVn (X)=\dim_{\CC}\mathscr{C}_0.
\]
\end{lemma}

\begin{proof}
Let $DI(0)\neq 0$. By Lemma \ref{lem1} we have
\[
\IVn (X)=\dim_{\CC}\mathscr{B}^I_0.
\]
If we map a class $[h]\in\mathscr{C}_0$ to the class $[h]'$ of $h$ in
$\mathscr{B}^I_0$ we obtain an isomorphism of $\CC$-algebras, since
by Lemma \ref{lem22} and Local Duality the following holds.
\begin{equation*}
\begin{split}
& [h]=0\;\;\textnormal{ in }\mathscr{C}_0\\
\Longleftrightarrow & hDF\in\OCn (f_1,\dots ,f_q,X_1,\dots ,X_{n-q})\\
\Longleftrightarrow & \RCn^X(ghDF)=0\;\;\textnormal{ for all }g\in\OCn\\
\Longleftrightarrow & \RCn\begin{bmatrix}
ghDI\\ X_{j_1}\dots X_{j_{n-q}}f_1\dots f_q\end{bmatrix}=0\;\;\textnormal{ for all }g\in\OCn\\
\Longleftrightarrow & hDI\in\OCn (X_{j_1},\dots X_{j_{n-q}},f_1,\dots ,f_q)\\
\Longleftrightarrow & [h]'=0\;\;\textnormal{ in }\mathscr{B}^I_0\\
\end{split}
\end{equation*}
\end{proof}

\begin{lemma}
\label{lem24}
Let $DI(0)\neq 0$. Then in $\mathscr{C}_0$ the equation
$c_{n-q}=\gamma_I$ holds.
\end{lemma}

\begin{proof}
By the transformation formula for residues and Lemma \ref{lem1} we have
\[
\RCn^X(DF\gamma_I )=\IVn (X)
\]
and further for $h\in\mathfrak{m}_{\OCn}$
\[
\RCn^X(DF\gamma_I h)=\RCn\begin{bmatrix}h\det
\left( \frac{\partial (X_{j_1},\dots ,X_{j_{n-q}},f_1,\dots ,f_q)}{\partial (z_1,\dots ,z_n)}\right) \\
X_{j_1}\dots X_{j_{n-q}}f_1\dots f_q\end{bmatrix}=0
\]
and therefore $DF\gamma_I$ generates the one-dimensional socle of $\mathscr{B}_0$. By the remarks of the introduction we have
\[
\RCn^X(DFc_{n-q})=\IVn (X),
\]
and therefore $DFc_{n-q}\neq 0$ in $\mathscr{B}_0$. Since $V$ is smooth there is a small deformation $X_t$
of $X$ tangent to $V$, so that $X_t$ has only simple zeros $p_i$ for sufficiently small $t$
on $V$ in a small neighbourhood of the origin. We can also assume that for these zeros
$DI(p_i)\neq 0$ holds. For $h\in\mathfrak{m}_{\OCn}$ we have
\begin{equation*}
\begin{split}
\RCn^X(DFc_{n-q}h)&= \RCn\begin{bmatrix} 
\sign\sigma_I DIc_{n-q}h\\X_{j_1}\dots X_{j_{n-q}}f_1\dots f_q\end{bmatrix}\\
&= \lim_{t\to 0}\sum_ih(p_i)\ind_{V,p_i}(X_t)\\
&= 0.
\end{split}
\end{equation*}
This follows from the continouos principle for residues, from the fact that the algebras
\[
\frac{\mathscr{O}_{\CC^n,p_i}}{(X_{t,j_1},\dots ,X_{t,j_{n-q}},f_1\dots ,f_q)}
\]
are one-dimensional and $\sign\sigma_I DIc_{n-q}(t)$ is a unit in these algebras by
the transformation formula. Therefore $c_{n-q}DF$ generates the one dimensional socle of
$\mathscr{B}_0$ too and since we have
\[
\RCn^X(DFc_{n-q})=\RCn^X(DF\gamma_I )
\]
it follows
$DF(c_{n-q}-\gamma_I )=0$ in $\mathscr{B}_0$
and therefore $c_{n-q}=\gamma_I$ in $\mathscr{C}_0$.
\end{proof}

\subsection{The real analytic situation}
\begin{lemma}
\label{lem6}
Let $(V^{\RR},0)$ be smooth and $l\colon\mathscr{C}^{\RR}_0\to\RR$ a linear form with
$l(c_{n-q})>0$. Then we have $\IVRn (X^{\RR})=\signa <,>_l$.
\end{lemma}

\begin{proof}
Let $DI (0)\neq 0$. With the implicit mapping theorem it is not hard to show that
the vector field $X^{\RR}|_{V^{\RR}}$ corresponds to
\[
Y:=X_{j_1}^{\RR}\circ\psi\frac{\partial}{\partial y_1}+\dots +
X_{j_{n-q}}^{\RR}\circ\psi\frac{\partial}{\partial y_{n-q}},
\]
where $\psi\colon (\RR^{n-q},0)\to (V^{\RR},0)$ is a diffeomorphism with
$\psi_{j_k}(y)=y_{j_k}$ for $k=1,\dots ,n-q$. By the chain rule one has
\[
\frac{\partial (\psi_{i_1},\dots ,\psi_{i_q})}{\partial (y_1,\dots ,y_{n-q})}=
-\left( \frac{\partial (f^{\RR}_1,\dots ,f^{\RR}_q)}{\partial (x_{i_1},\dots ,x_{i_q})}\circ\psi\right)^{-1}
\frac{\partial (f^{\RR}_1,\dots ,f^{\RR}_q)}{\partial (x_{j_1},\dots ,x_{j_{n-q}})}\circ\psi .
\]
Applying the chain rule again we get
\begin{multline*}
\frac{\partial (Y_1,\dots ,Y_{n-q})}{\partial (y_1,\dots ,y_{n-q})}=
\frac{\partial (X^{\RR}_{j_1},\dots ,X^{\RR}_{j_{n-q}})}{\partial (x_{j_1},\dots ,x_{j_{n-q}})}\circ\psi\\
- \frac{\partial (X^{\RR}_{j_1},\dots ,X^{\RR}_{j_{n-q}})}{\partial (x_{i_1},\dots ,x_{i_q})}\circ\psi
\left( \frac{\partial (f^{\RR}_1,\dots ,f^{\RR}_q)}{\partial (x_{i_1},\dots ,x_{i_q})}\circ\psi\right)^{-1}
\frac{\partial (f^{\RR}_1,\dots ,f^{\RR}_q)}{\partial (x_{j_1},\dots ,x_{j_{n-q}})}\circ\psi .
\end{multline*}
A well known Lemma from linear algebra states
\[
\det\begin{pmatrix}A & B\\ C & D\end{pmatrix}=\det A\cdot\det (D-CA^{-1}B),
\]
where $A$ and $D$ are squared and $A$ is invertible.
The application of this Lemma shows that the determinant of the Jacobian of $Y$ is given by
$\gamma_I\circ\psi$. By the Eisenbud-Levine Theorem it follows that for any linear form
\[
\varphi\colon\frac{\mathscr{E}_{\RR^{n-q},0}}{(Y_1,\dots ,Y_{n-q})}\to\RR
\]
with $\varphi (\gamma_I\circ\psi )>0$ the statement $\IVRn (X)=\signa <,>_{\varphi}$ holds. The isomorphism of algebras
given by $\psi$ shows that we have for any linear form
$\Phi\colon\mathscr{B}_0^{I\RR}\to\RR$ with $\Phi (\gamma_I )>0$
the formula $\IVRn (X)=\signa <,>_{\Phi}$
and this is also true in $\mathscr{C}_0^{\RR}$, because the isomorphism of algebras
in Lemma \ref{lem3} also gives an isomorphism of the corresponding real algebras.
On the other hand one has in $\mathscr{C}_0$ the equation $c_{n-q}=\gamma_I$ and this equation also holds in
$\mathscr{C}_0^{\RR}$. Therefore the statement follows.
\end{proof}

\section{An algebraic formula for the complex index}

\subsection{Good vector fields}
First we want to prove a sufficient criterion for good vector fields.

\begin{satz}[(Sufficient criterion for good vector fields)]
Let all coefficients of the matrix $C$ be contained in the ideal generated by the maximal
minors of the Jacobian of $f$ in $\OCn$. Then $X$ is a good vector field.
\end{satz}

\begin{proof}
We prove a bit more: There is a deformation $X_t$ of $X$ so that
$X_t(f-t)=C(f-t)$ holds. For $(i_1,\dots ,i_q)\in\{1,\dots ,n\}^q$ and
$(j_1,\dots ,j_{q-1})\in\{1,\dots ,n\}^{q-1}$ define
\[
f_{i_1,\dots ,i_q}:=\det\left(\frac{\partial (f_1,\dots ,f_q)}{\partial (z_{i_1},\dots ,z_{i_q})}\right)
\text{ and }
f^k_{j_1,\dots ,j_{q-1}}:=\det\left(\frac{\partial (f_1,\dots ,\hat{f_k},\dots ,f_q)}{\partial (z_{j_1},\dots ,z_{j_{q-1}})}\right) .
\]
If $C=(c_{l,m})$ let
\[
c_{l,m}=\sum_{(i_1,\dots ,i_q)\in\{1,\dots ,n\}^q}c_{i_1,\dots ,i_q}^{l,m}f_{i_1,\dots ,i_q}.
\]
We set for $k=1,\dots ,q-1$
\[
\delta_{j_1,\dots ,j_{q-1}}^{i,j}:= \delta_{j_1,\dots ,j_{k-1},i,j_{k+1},\dots ,j_{q-1}}^{j_k,j}:=
\delta_{i, j_1,\dots ,j_{q-1}}^j
\]
and further
\[
\delta_{i_1,\dots ,i_q}^l:=-\sum_{m=1}^q\sum_{(i_1,\dots ,i_q)\in\{1,\dots ,n\}^q}
c_{i_1,\dots ,i_q}^{l,m}t_m.
\]
We define the deformation by
\[
X_{t,i}:=X_i+\sum_{j=1}^q\sum_{(j_1,\dots ,j_{q-1})\in\{1,\dots ,n\}^{q-1}}
\delta_{j_1,\dots ,j_{q-1}}^{i,j}(-1)^{j+1}f_{j_1,\dots ,j_{q-1}}^j.
\]
Then we have
\begin{equation*}
\begin{split}
\sum_{i=1}^n\frac{\partial f_l}{\partial z_i}X_{t,i}&=
\sum_{m=1}^nc_{l,m}f_m+\sum_{(i_1,\dots ,i_q)\in\{1,\dots ,n\}^q}\delta_{i_1,\dots ,i_q}^lf_{i_1,\dots i_q}\\ &=
\sum_{m=1}^nc_{l,m}(f_m-t_m).
\end{split}
\end{equation*}
\end{proof}

In the case of a hypersurface this means if $c$
is contained in the ideal generated by the partials of $f$
the vector field is good. If
\[
c=\alpha_1\frac{\partial f}{\partial z_1}+\dots +\alpha_n\frac{\partial f}{\partial z_n}
\]
the deformation is simply defined as
$X_{t,i}:=X_i-t\alpha_i$. We have
\[
X_t(f-t)=cf-\sum_{i=1}^nt\alpha_i\frac{\partial f}{\partial z_i}=c(f-t).
\]

\subsection{The socle of $\mathscr{C}_0$}
\begin{lemma}
\label{lem5}
$\RVn^X(\cdot )$ defines a linear form on $\mathscr{C}_0$ so that
\[
\RVn^X(hg)=0\;\textnormal{ for all }\; h\in\OCn\Rightarrow  g=0 \textnormal{ in }\mathscr{C}_0
\]
holds.
\end{lemma}

\begin{proof}
We have $\RVn^X (h)=\RCn^X(hDF)$ and this means that the residue
$\RVn^X (\cdot )$ vanishes on $\ann_{\mathscr{B}_0}(DF)$. Furthermore $\RVn^X(hg)=0$ for all $h\in\OCn$ implies
$\RCn^X(hgDF)=0$ for all $h\in\OCn$. Now Local Duality gives
\[
gDF\in\OCn (X_1,\dots ,X_{n-q},f_1,\dots ,f_q)
\]
and therefore $g\in\ann_{\mathscr{B}_0}(DF)$.
\end{proof}

\begin{satz}
The socle of $\mathscr{C}_0$ is generated by the class of $c_{n-q}$.
\end{satz}

\begin{proof}
Lemma \ref{lem5} states that $\RVn^X(\cdot )$ induces a non-degenerate bilinear form on $\mathscr{C}_0$ which means
that $\mathscr{C}_0$ has an one-dimensional socle if $\mathscr{C}_0$ is not trivial.
On the other hand we have for any
$h\in\mathfrak{m}_{\OCn}$
\[
\RVn^X(hc_{n-q})=\lim_{t\to 0}\sum_ih(p_i)\ind_{V_t,p_i}(X_t)=0,
\]
which means that $c_{n-q}$ generates the socle.
\end{proof}

\subsection{Proof of Theorem \ref{them1}}
For a good vector field $\IVn (X)$ is the sum of Poincar\'{e}-Hopf indices of $X_t$
on a smooth fibre $V_t$ of $f$, where
one sums over all zeros of $X_t$ which tend to zero. This follows directly from the definition of the
index.

If $q=n-1$ we denote by $m_i$ the minor of the Jacobian matrix of $f$ obtained by cancelling the $i$-th column. We have $DF=m_1$
of course. Recall that we always assume the coordinates to be chosen so that $(f_1,\dots , f_{n-1},X_1)$ is a regular
$\OCn$-sequence. For simplicity we also may assume that $m_i(0)=0$ for $i=1,\dots ,n$.

\begin{lemma}
\label{coor}
Let $q=n-1$. Then\\
(i) $(f_1,\dots ,f_{n-1},DF)$ is a regular $\OCn$-sequence.\\
(ii) The dimension of $\mathscr{C}_0$ does not depend on the choice of coordinates, provided that
the $\OCn$-sequence $(f_1,\dots ,f_{n-1},X_1)$ is regular in each system of coordinates.\\
(iii) After linear changes of coordinates in $\CC^n$ and $\CC^{n-1}$ one can assume that
$(f_1,\dots ,f_{n-2},DF,m_2)$ is a regular $\OCn$-sequence.
\end{lemma}

\begin{proof}
(i) All computations are done in the ring
$\mathscr{O}_{V,0}=\OCn /(f_1,\dots ,f_{n-1})$. We have to show that $DF$ is not
a zero divisor in this ring. Applying Cramer's rule to the equation $Xf=Cf$ we obtain the equations
\begin{equation}
\label{eqcr}
(-1)^im_iX_1=-DFX_i
\end{equation}
for $i=1,\dots ,n$. Now let $gDF=0$ in $\mathscr{O}_{V,0}$. Multiplication with $X_i$, equation \ref{eqcr} and using that
$X_1$ is not a zero divisor in $\mathscr{O}_{V,0}$ gives $gm_i=0$ for all $i=1,\dots ,n$. The ring obtained by dividing
$\mathscr{O}_{V,0}$ by all $m_i$ is artinian and therefore there must be complex numbers $\alpha_1 ,\dots ,\alpha_n$ so that
$h:=\alpha_1 m_1+\dots +\alpha_n m_n$ is not a zero divisor in $\mathscr{O}_{V,0}$. On the other hand we have $gh=0$ in
$\mathscr{O}_{V,0}$ and therefore $g=0$ in $\mathscr{O}_{V,0}$.\\
(ii) Let $\phi\colon (\CC^n,0)\to (\CC^n,0)$, $\phi(y)=z$, be biholomorphic and $\psi :=\phi^{-1}$. We denote by $Y$ the vector
field computed in $y$-coordinates and by $DF^y$ the minor computed in $y$-coordinates.
Standard computations give
\[
\begin{pmatrix} Y_1\circ\psi \\ \vdots \\ Y_n\circ\psi \end{pmatrix}
=\frac{\partial (\psi_1,\dots ,\psi_n )}{\partial (z_1,\dots ,z_n)}
\begin{pmatrix} X_1\\ \vdots \\ X_n \end{pmatrix}
\]
and
\[
DF^y\circ\psi = \sum_{j=1}^n(-1)^{j+1}(\det D\phi )\circ\psi\;\frac{\partial \psi_1}{\partial z_j}m_j.
\]
Set $\mathscr{B}'_0:=\OCn /(f_1,\dots ,f_{n-1},Y_1\circ\psi )$ and $\mathscr{C}'_0:=\mathscr{B}'_0/\ann_{\mathscr{B}'_0}
(DF^y\circ\psi )$.
We construct an epimorphism
$\varphi\colon\mathscr{C}_0\to\mathscr{C}'_0$.
For any $g\in\OCn$ we also denote by $g$ the classes of $g$ in these algebras and define $\varphi (g):=g$. Again all
computations are done in the ring $\mathscr{O}_{V,0}$. We want to show that $\varphi$ is well defined. Let $gDF=\alpha X_1$. Multiplication
with $m_i$, usage of equation \ref{eqcr} and usage of the fact that $DF$ is not a zero divisor in $\mathscr{O}_{V,0}$ we obtain
$gm_i=(-1)^{i+1}\alpha X_i$ for $i=1,\dots ,n$. Then for each $i$ we get
\begin{equation*}
\begin{split}
gDF^y\circ\psi &= (\det D\phi )\circ\psi\;\sum_{i=1}^n(-1)^{i+1}\frac{\partial \psi_1}{\partial z_i}gm_i\\
&= (\det D\phi )\circ\psi\;\sum_{i=1}^n\alpha X_i \frac{\partial \psi_1}{\partial z_i}\\
&= \alpha (\det D\phi )\circ\psi\; Y_1\circ\psi .
\end{split}
\end{equation*}
This shows that $\varphi$ is well defined. The surjectivity is obvious and the other direction analogous. Therefore
$\dim_{\CC}\mathscr{C}_0$ does not depend on the choice of coordinates in $\CC^n$. If one considers changes of coordinates
in the image space $\CC^{n-1}$ the invariance of $\dim_{\CC}\mathscr{C}_0$ is obvious.\\
(iii) After a general linear change of coordinates in $\CC^{n-1}$, see \cite{l}, one can assume that $(f_1,\dots ,f_{n-2})$ defines
an icis and the $1$-form $df_{n-1}$ has an isolated zero on this icis. Now Lemma 3.4 in \cite{k} shows that after a linear
change of coordinates in $\CC^n$ $(f_1,\dots ,f_{n-2}, DF,m_2)$ is a regular $\OCn$-sequence.
\end{proof}

\begin{proof}[Proof of Theorem \ref{them1}]
Let $q=n-1$ and consider good deformations $X_t$ of $X$: For small neighbourhoods $U$ resp. $T$ of the origins in $\CC^n$ resp. $\CC^{n-1}$
consider $Z\subset U\times T$ defined by
\[
Z:=\{ f_{t,1}=\dots =f_{t,n-1}= X_{t,1}=0\},
\]
where we have set $f_{t,i}:=f_i-t_i$. Let $\pi\colon Z\to T$ be the finite projection. We define
\[
\mathscr{B}_{t,p}:=\frac{\mathscr{O}_{\CC^n,p}}{(f_{t,1},\dots ,f_{t,n-1},X_{t,1})},\;\;
\mathscr{B}_t:=\bigoplus_{p\in\pi^{-1}(t)}\mathscr{B}_{t,p},\;\;
\mathscr{B}:=\bigcup_{t\in T}\mathscr{B}_t.
\]
$\mathscr{B}$ has the natural structure of a holomorphic vector bundle over $T$ which is induced by the locally free sheaf
$\pi_{*}\mathscr{O}_Z$. Similarly we define
\[
\mathscr{C}_{t,p}:=\frac{\mathscr{B}_{t,p}}{\ann_{\mathscr{B}_{t,p}}(DF)},\;\;
\mathscr{C}_t:=\bigoplus_{p\in\pi^{-1}(t)}\mathscr{C}_{t,p},\;\;\mathscr{C}:=\bigcup_{t\in T}\mathscr{C}_t.
\]
We want to show that $\mathscr{C}$ has the natural structure of a holomorphic vector bundle of rank
$\IVn (X)$ over $T$.

By Lemma \ref{coor} (ii) we may assume that the coordinates are chosen as in Lemma $\ref{coor}$ (iii). Via $\pi$ we view
$\mathscr{O}_{Z,0}/(DF)$ as finitely generated $\mathscr{O}_{T,0}$-module and claim that
\[
\dep_{\mathscr{O}_{T,0}}\frac{\mathscr{O}_{Z,0}}{(DF)}=n-1.
\]
For $k=1,\dots ,n-1$ let
\[
t_kg=0\;\;\textnormal{in }\frac{\mathscr{O}_{Z,0}}{(DF,t_1,\dots ,t_{k-1})}.
\]
This means that there are representatives with
\[
t_kg=\alpha X_{t,1}\;\;\textnormal{in }\frac{\mathscr{O}_{\CC^{2n-1},0}}{(f_{t,1},\dots ,f_{t,n-1},t_1,\dots ,t_{k-1},DF)}.
\]
Applying Cramer's rule to the tangency equation $Xf=Cf$ we obtain
\[
t_kgm_2=0\;\;\textnormal{in }\frac{\mathscr{O}_{\CC^{2n-1},0}}{(f_{t,1},\dots ,f_{t,n-1},t_1,\dots ,t_{k-1},DF)}.
\]
By Lemma \ref{coor} (i) and the choice of coordinates $t_k$ and $m_2$ are no zero divisors in the last algebra since
$(f_{t,1},\dots ,f_{t,n-1},t_1,\dots ,t_{n-2},DF,m_2)$ is a regular $\mathscr{O}_{\CC^{2n-1},0}$-sequence.
This shows the claim.

Now the Syzygy Theorem and the Auslander-Buchsbaum formula show that $\mathscr{O}_{Z,0}/(DF)$ is a
free $\mathscr{O}_{T,0}$-module. We have an exact sequence of $\mathscr{O}_{T,0}$-modules
\[
0\to\frac{\mathscr{O}_{Z,0}}{\ann_{\mathscr{O}_{Z,0}}(DF)}\to \mathscr{O}_{Z,0}\to
\frac{\mathscr{O}_{Z,0}}{(DF)}\to 0.
\]
The Depth Lemma, see \cite[6.5.18]{jp}, shows that
\[
\dep_{\mathscr{O}_{T,0}}\frac{\mathscr{O}_{Z,0}}{\ann_{\mathscr{O}_{Z,0}}(DF)}=n-1,
\]
which means by the Syzygy Theorem and the Auslander-Buchsbaum formula again that
the $\mathscr{O}_{T,0}$-module
$\mathscr{O}_{Z,0}/\ann_{\mathscr{O}_{Z,0}}(DF)$ is free.
We have seen that the coherent $\mathscr{O}_T$-module
\[
\mathscr{F}:=\pi_{*}\mathscr{O}_Z/\ann_{\mathscr{O}_Z}(DF)
\]
is free for $T$ chosen small enough. By \cite{d} this is equivalent to the statement that the function
$\nu\colon T\to\mathds{N}$ defined by
\[
\nu (t):=\dim_{\CC}\mathscr{F}_t\otimes_{\mathscr{O}_{T,t}}\CC
\]
is constant. Now from the exact sequence
\[
0\to \mathscr{F}_t\to (\pi_{*}\mathscr{O}_Z)_t\to (\pi_{*}\mathscr{O}_Z/(DF))_t\to 0
\]
we obtain by tensoring with $\CC$ the last part of the long exact sequence of torsion
\[
0\to\mathscr{F}_t\otimes_{\mathscr{O}_{T,t}}\CC\to\mathscr{B}_t\to \mathscr{B}_t/(DF)\to 0,
\]
where we have used that $(\pi_{*}\mathscr{O}_Z/(DF))_t$ is a free $\mathscr{O}_{T,t}$-module which
is equivalent, see \cite{d}, to the statement that
\[
\tor^{\mathscr{O}_{T,t}}_1((\pi_{*}\mathscr{O}_Z/(DF))_t,\CC )=0.
\]
Since we also have an exact sequence
\[
0\to\mathscr{C}_t\to\mathscr{B}_t\to\mathscr{B}_t/(DF)\to 0
\]
this means that $\nu (t)=\dim_{\CC}\mathscr{C}_t$ for all $t\in T$. By Lemma
\ref{lem3} $\nu (t)$ equals for regular $t$ the sum of Poincar\'{e}-Hopf indices of $X_t$ on
the $t$-fibre of $f$ for $t\neq 0$ which is equal to $\IVn (X)$ and therefore
$\dim_{\CC}\mathscr{C}_0=\IVn (X)$. The map $\cdot DF\colon\mathscr{B}\to\mathscr{B}$
between vector bundles has constant rank and provides $\mathscr{C}$ with the natural
structure of a holomorphic vector bundle.
\end{proof}

\subsection{Examples and Remarks}
We want to give a few examples of good vector fields on curves.
We can always take the exterior pruduct of the rows of the Jacobian matrix
of $f$. In this way one obtains a vector field with isolated zero on the singularity, where the matrix $C$ is
trivial and the index equals zero.

An example of a family of non-trivial good vector fields on a plane curve is the following:
\[
D_k:\;f=x^2y+y^{k-1},\;k\geq 4
\]
and
\[
X=(k-2)x^{m+1}\frac{\partial}{\partial x}+2x^my\frac{\partial}{\partial y},\;m\geq 3
\]
with $c=2(k-1)x^m$.
We have an exact sequence
\[
0\xrightarrow{} \ann_{\mathscr{B}_0}(DF)\xrightarrow{} \mathscr{B}_0 \xrightarrow{\cdot DF} \mathscr{B}_0 \xrightarrow{} \frac{\mathscr{B}_0}{(DF)}\xrightarrow{} 0
\]
and therefore
\[
\dim_{\CC}\ann_{\mathscr{B}_0}(DF)=\dim_{\CC}\frac{\mathscr{B}_0}{(DF)}.
\]
It is easy to compute that $\dim_{\CC}\mathscr{B}_0/(DF)=2(k-1)$ and
$\dim_{\CC}\mathscr{B}_0=(k-1)(m+1)$ holds and therefore $\IVn (X)=(k-1)(m-1)$.

An example of a family of good vector fields on a space curve is
\[
f_1:=x^2+y^2+z^2\;\; ,f_2:=xy
\]
with
\[
X:=z^l(x-y)\left( x\frac{\partial}{\partial x}+y\frac{\partial}{\partial y}+
z\frac{\partial}{\partial z}\right)\;\; ,l\geq 1 .
\]
Here we have
\[
C=
\begin{pmatrix}
2z^l(x-y) & 0\\
0 & 2z^l(x-y)
\end{pmatrix}.
\]
In this case $(f_1, f_2, X_3)$ is a regular sequence and we have to permute the coordinates. This means
that the index is given by the dimension of the factor space obtained from factoring the algebra
$\mathscr{O}_{\CC^3,0}/(f_1,f_2, X_3)$ by the annihilator of
\[
\det
\begin{pmatrix}
2x & 2y\\
y & x
\end{pmatrix}.
\]
In this case we get a non trivial index, which one may to compute with computer algeba programs
such as Singular \cite{gps}. Note that $\dim_{\CC}\mathscr{C}_0$ is always equal to
$\dim_{\CC}\mathscr{B}_0-\dim_{\CC}\mathscr{B}_0/(DF)$. The last two dimensions are easily to compute
with computer algebra.

We now want to explain why Theorem \ref{them1} is false in the general case of an complete intersection. Consider the case
$n=3$ and $q=1$ and set $f_i:=\partial f/\partial z_i$ for $i=1,2,3$. We use the notations as in the proof of Theorem
\ref{them1}, i.e. $\mathscr{O}_{Z,0}:=\mathscr{O}_{4,0}/(f-t, X_{t,1}, X_{t,2})$ and $(f,X_1, X_2)$ is a regular
sequence in $\mathscr{O}_{3,0}$, and consider the coherent
$\mathscr{O}_T$-module $\pi_{*}\mathscr{O}_Z/(f_3)$ where $T$ is a small
neighbourhood of $0$ in $\CC$. After shrinking $T$ we can assume that this sheaf is locally free over all $t\neq 0$.
Obviously $t$ is not a zero divisor in $\mathscr{F}_0$ and therefore we have a law of conservation of numbers as in the proof of
\ref{them1} and we get
\[
\IVn (X)=\dim_{\CC}\mathscr{F}_0\otimes_{\mathscr{O}_{T,0}}\CC.
\]
We further have the exact sequence
\[
0\to\tor_1^{\mathscr{O}_{T,0}}(\pi_{*}\mathscr{O}_{Z,0}/(f_3),\CC )\to\mathscr{F}_0\otimes_{\mathscr{O}_{T,0}}\CC\to
\mathscr{B}_0\to\mathscr{B}_0/(f_3)\to 0
\]
and this means
\[
\IVn (X)=\dim_{\CC}\mathscr{C}_0 + \dim_{\CC}\tor_1^{\mathscr{O}_{T,0}}(\pi_{*}\mathscr{O}_{Z,0}/(f_3),\CC ).
\]
We now show that $\dim_{\CC}\tor_1^{\mathscr{O}_{T,0}}(\pi_{*}\mathscr{O}_{Z,0}/(f_3),\CC )>0$ if the hypersurface is not smooth
which is equivalent to the statement that $t$ is a zero divisor in $\mathscr{O}_{Z,0}/(f_3)$. Set
$\mathscr{O}_{4,0}:=\CC\{ z_1,z_2,z_3,t\}$. The functions $(f-t, f_1,f_2,f_3)$ define an isolated zero in $(\CC^4,0)$ and therefore these
define a regular $\mathscr{O}_{4,0}$-sequence. We have
\[
f_1X_{t,1}+f_2X_{t,2} =0
\]
in the ring $\mathscr{O}_{4,0}/(f-t,f_3)$. It follows that there are $\gamma_1,\gamma_2\in\mathscr{O}_{4,0}$ such that
$X_{t,1}=\gamma_1f_2$ and $X_{t,2}=\gamma_2f_1$ in $\mathscr{O}_{4,0}/(f-t,f_3)$ holds. Inserting this in the tangency equation
shows that there is a $\gamma\in\mathscr{O}_{4,0}$ so that
$X_{t,1}=\gamma f_2$ and $X_{t,2}=-\gamma f_1$ in $\mathscr{O}_{4,0}/(f-t,f_3)$ holds. On the other hand
$\mathscr{O}_{Z,0}/(f_3,t)$ is artianian and therefore $(f-t,\gamma ,f_3, t)$ is a weak regular $\mathscr{O}_{4,0}$-sequence.
This means that we have an exact sequence
\[
0\rightarrow\frac{\mathscr{O}_{4,0}}{(f-t,f_1,f_2,f_3)}\stackrel{\cdot\gamma}{\rightarrow}\frac{\mathscr{O}_{4,0}}{(f-t,\gamma f_1,\gamma f_2, f_3)}
\rightarrow\frac{\mathscr{O}_{4,0}}{(f-t,\gamma ,f_3)}\rightarrow 0
\]
with
\[
\frac{\mathscr{O}_{4,0}}{(f-t,\gamma f_1, \gamma f_2, f_3)}\cong\frac{\mathscr{O}_{Z,0}}{(f_3)}.
\]
The sequence shows that $\dep_{\mathscr{O}_{T,0}}\mathscr{O}_{Z,0}/(f_3)=0$:
If $\gamma (0)\neq 0$ this is obviously. If $\gamma (0)=0$ we can apply the Depth Lemma using
\[
\dep_{\mathscr{O}_{T,0}}\frac{\mathscr{O}_{4,0}}{(f-t, f_1, f_2, f_3)}=0
\]
and
\[
\dep_{\mathscr{O}_{T,0}}\frac{\mathscr{O}_{4,0}}{(f-t, \gamma ,f_3)}=1.
\]

\section{A signature formula for the real index}
Now let
\[
(V^{\RR},0):=(\{ f^{\RR}_ 1=\dots =f^{\RR}_ q =0\} ,0)\subset (\RR^n,0)
\]
be a geometric complete intersection of dimension $n-q$, denote by $V$ and $f$ the complexifications and assume
that $f$ defines an ICIS.
Furthermore let the real analytic vector field $X^{\RR}$ be tangent to $(V_{\RR},0)$ with an algebraic
isolated zero on $(V^{\RR},0)$. 
As before $T$ is a small neighbourhood
of the origin in $\CC^q$ and let $T^{\RR}$ be the corresponding subset in $\RR^q$. We also assume that
$X^{\RR}$ is good in the real analytic sense. We keep the notations of the previous section
for all complexifications and denote for real $t$ the real algebra
corresponding to $\mathscr{C}_{t,p}$ by $\mathscr{C}_{t,p}^{\RR}$, if $X^{\RR}|_{V_t^{\RR}}(p)=0$ and $p\in\RR^n$ holds.

\subsection{Proof of Theorem \ref{them2}}
First we prove a law of conservation of numbers for the signature. Theorem \ref{them2}
follows as a corollary then. We also remark that the sufficient criterion for good vector fields also holds in the real
analytic case. The proof is word for word as in the complex case.

\begin{satz}
\label{prop22}
Let $q=n-1$, $X^{\RR}$ be a good vector field and $l\colon\mathscr{C}_0^{\RR}\to\RR$ a linear form
with $l(c_1)>0$ and for any regular value
$t\in T^{\RR}$ of $f^{\RR}$ and any $p$ with $X_t^{\RR}|_{V_t^{\RR}}(p)=0$ let $l_{t,p}\colon
\mathscr{C}_{t,p}^{\RR}\to\RR$ be a linear form with $l_{t,p}(c_{t,p,1})>0$. Then
\[
\signa <,>_l =\sum_{\{X_t^{\RR}|_{V_t^{\RR}}(p)=0\}}\signa <,>_{l_{t,p}},
\]
where the sum goes over the zeros tending to zero.
\end{satz}

\begin{proof}
We consider the vector bundle $\mathscr{C}$ over $T$ defined in the proof of \ref{them1} and denote the map
given by complex conjugation by $\tau$.
For any $t\in T^{\RR}$ we consider the set of invariant multigerms $h\in\mathscr{C}_t$. These are the multigerms with
$\tau\circ h=h\circ\tau$. We denote this set by $\mathscr{C}_t^{\RR}$. We have
\begin{equation}
\label{equa1}
\mathscr{C}_t^{\RR}=\left( \oplus_kD_k\right) \oplus\left( \oplus_l E_l\right),
\end{equation}
where each component $D_k$ corresponds to an algebra $\mathscr{C}^{\RR}_{t,p_k}$ for a real zero $p_k$ of $X|_V$ and where
each component
\[
E_l=(\mathscr{C}_{t,q_l}\oplus\mathscr{C}_{t,\overline{q_l}})^{\RR}
\]
corresponds to a pair of complex conjugated zeros and $(\mathscr{C}_{t,q_l}\oplus\mathscr{C}_{t,\overline{q_l}})^{\RR}$ is the subset
of invariant elements of $(\mathscr{C}_{t,q_l}\oplus\mathscr{C}_{t,\overline{q_l}})$. It consists of elements of the form
\[
h=\sum a_Iz^I + \sum\overline{a_I}z^I.
\]
Here $q_l$ resp. $\overline{q_l}$ are not real of course. If $\mu$ is the real dimension of
$\mathscr{C}_t^{\RR}$ then $\mu$ is given by $\dim_{\CC}\mathscr{C}_0$. The set
$\mathscr{C}^{\RR}:=\cup_{t\in T^{\RR}}\mathscr{C}_t^{\RR}$ has for $T$ chosen small enough
the natural structure of a real analytic vector bundle of rank $\mu$ over $T^{\RR}$. We can continue
$l$ real analytically to a family $l_t$ and obtain a real analytic family
of non-degenerate bilinear forms $<,>_{l_t}$. Equation \ref{equa1} gives an orthogonal decomposition.
By dividing the algebra $E_l$ by its maximal ideal one obtains $\CC$ and therefore
$E_l$ contributes nothing to the signature, see \cite{el}. Therefore the signature of
$<,>_{l_t}$ is the sum of signatures of $<,>_{l_t,p}$ that are defined as the restrictions to the components
$D_k$. On the other hand we have $l_{t,p}(c_{t,p,1})>0$ and therefore the claim follows
by continuity of signatures and by the theorem of Eisenbud and Levine if we choose a fixed regular value $t\in T^{\RR}$ of $f$.
\end{proof}

\begin{proof}[Proof of Theorem \ref{them2}]
For a good vector field the index counts the sum of indices of a good deformation of the
vector field on a regular fibre in a neighbourhood of the origin by the properties of
the real index given in Theorem 2.10 in \cite{asv}. Now the claim follows from Lemma
\ref{lem6} and Proposition \ref{prop22}.
\end{proof} 

We want to consider an example.
Let $f^{\RR}(x,y):=x^2-y^2$ and $X^{\RR}:=x^2\partial /\partial x+ xy\partial /\partial y$.
A good deformation is given by
\[
X_t^{\RR}:=(x^2-t)\frac{\partial}{\partial x}+xy\frac{\partial}{\partial y}
\]
with $c_t=c=2x$. Set $F_t:=V^{\RR}_t\cap\overline{B_{\delta}}$ where $B_{\delta}$ is a small
ball around the origin in $\RR$. $F_t$ consists of two branches of a hyperbola and we have $\chi (F_t)=2$. If $l$
is a linear form as in Theorem \ref{them2} we obtain as easily to see
$\signa <,>_l=0$. Let $t=1$,
$B_{\delta}:=\{x^2+y^2=3\}$ and $F:=F_1$. $X^{\RR}$ deforms to
\[
\tilde{X}^{\RR}:=(x^2-1)\frac{\partial}{\partial x}+xy\frac{\partial}{\partial y}.
\]
then. The boundary points of $F$ are $P_1=(\sqrt{2},1),\;P_2=(\sqrt{2},-1),\;
P_3=(-\sqrt{2},-1)$ and $P_4=(-\sqrt{2},1)$. At the points $P_1$ und $P_2$ the vector field $\tilde{X}^{\RR}$ points
outwards but inwards at the points $P_3$ and $P_4$. From the symmetry
of the problem (only the directions of $\tilde{X}^{\RR}$ are not symmetric) we find that the sum
of the indices of $\tilde{X}^{\RR}$ vanishes on $F$ and this is what Theorem \ref{them2} says. This can also
be computed explicitely: The zeros of $\tilde{X}^{\RR}$ on $F$ are
$(-1,0)$ und $(1,0)$. We can parametrize both branches via $\varphi_{\pm}(s):=(\pm\sqrt{1+s^2},s)$
and write $\tilde{X}^{\RR}$ in the coordinate given by $s$. One immediately sees that
the index in $(-1,0)$ has the value $-1$ and the value $1$ in $(1,0)$.

If we want to count the Euler characteristic of $F_t$ we have to choose a good vector field whose deformation
points outwards at all boundary points. This means that we have to choose a good vector field which
points outwards at all boundary points of the intersection of the singular fibre with a small
closed ball around the origin.

\subsection{Relations to results of G\'{o}mez-Mont and Mardesi\'{c}}

G\'{o}mez-Mont and Mardesi\'{c} have proven similar signature formulas
in the articles \cite{gm} and \cite{gm2}.
These formulas hold for vector fields on isolated hypersurface singularities which have not
only an isolated zero on the variety but also in the ambient space. We want to compare these
formulas with our formula for vector fields on plane curves. 
Let $X$ be a real analytic vector field in $(\RR^n,0)$ with an isolated zero and $(V,0):=
(\{f=0\},0)\subset(\RR^n,0)$ an odd-dimensional hypersurface with algebraic isolated singularity. Further let
$X$ be tangent to $V$, i.e. $Xf=cf$. We omit the upper $\RR$ to indicate that we are
working in the real analytic category. Define
\[
\mathds{A}:=\frac{\ERn}{(f_1,\dots f_n)} \text{ and } \mathds{B}:=\frac{\ERn}{(X_1,\dots ,X_n)}.
\]
Here the $f_i$ are the partials of $f$. Let $H_f$ be the Hessian determinant of $f$. $\det DX$
and $H_f$ generate the socles of $\mathds{B}$ resp. $\mathds{A}$. Now we have
well determined classes
\[
H_f^{\text{rel}}:=\frac{H_f}{c}\in\frac{\mathds{A}}{\ann_{\mathds{A}}(c)},\;\;\;\;\;
\det DX^{\text{rel}}:=\frac{\det DX}{c}\in\frac{\mathds{B}}{\ann_{\mathds{B}}(c)}
\]
defined in the obvious way, which generate the socles of these algebras. Let
\[
l_1\colon \frac{\mathds{A}}{\ann_{\mathds{A}}(c)}\to\RR,\;\;\;\;\;
l_2\colon \frac{\mathds{B}}{\ann_{\mathds{B}}(c)}\to\RR
\]
be linear forms with $l_1 (H_f^{\text{rel}})>0$ and $l_2 (det DX^{\text{rel}})>0$. We have
\begin{thm}[G\'{o}mez-Mont, Mardesi\'{c}]
\label{gm}
\[
\IVn (X)=\signa <,>_{l_2}- \signa <,>_{l_1}.
\]
\end{thm}
To compare this result with our theorem we additionally assume $n=2$, the
$\mathscr{E}_{2,0}:=\RR\{x,y\}$-sequence $(f,X_1)$ to be regular and the vector field to
be good.  
We first give an explicit construction of all vector fields fulfilling both
conditions: 
Denote by $X_t$ the good deformation. We also set
$\mathscr{E}_{3,0}:=\RR\{x,y,t\}$,
$\mathscr{E}_{V,0}:=\mathscr{E}_{3,0}/(f-t)$. The tangency condition gives
\[
X_{t,1}f_1+X_{t,2}f_2=0 \text{ in } \mathscr{E}_{V,0}.
\]
Since $(f_1,f_2)$ is a regular
$\mathscr{E}_{V,0}$-sequence, it follows immediately that there are $\tilde{\gamma},
\tilde{\delta_1},\tilde{\delta_2}\in\mathscr{E}_{3,0}$ so that
\[
X_t=(\tilde{\gamma}f_2+\tilde{\delta_1} (f-t))\frac{\partial}{\partial x}+
(-\tilde{\gamma}f_1+\tilde{\delta_2}(f-t))\frac{\partial}{\partial y}.
\]
Setting $t=0$ we obtain that there are $\delta_1,\delta_2,\gamma\in\mathscr{E}_{2,0}$ so that
\[
X=(\gamma f_2+\delta_1 f))\frac{\partial}{\partial x}+
(-\gamma f_1+\delta_2 f)\frac{\partial}{\partial y}.
\]
We have $c=\delta_1f_1+\delta_2f_2$ and this means $\signa <,>_{l_1}=0$. We now claim that
$\ann_{\mathds{B}}(c)=\mathds{B}(\gamma ,f)$. Using
\begin{equation*}
\begin{split}
cf &= f_1X_1+f_2X_2\\
c\gamma &= \delta_2X_1 - \delta_1X_2
\end{split}
\end{equation*}
we obtain $\mathds{B}(\gamma ,f)\subset \ann_{\mathds{B}}(c)$. Now let $cg=\alpha_1X_1+\alpha_2X_2$.
Multiplication with $f$ gives
\[
(f_1g-\alpha_1f)X_1+(f_2g-\alpha_2f)X_2=0.
\]
Since $(X_1,X_2)$ is a regular $\mathscr{E}_{2,0}$-sequence there must be an $h\in\mathscr{E}_{2,0}$
with $f_2g-\alpha_2f=hX_1$. Since $(f,X_1)$ is a regular $\mathscr{E}_{2,0}$-sequence this
also holds for $f,f_2$. Therefore we have $g=h\gamma$ in $\mathscr{E}_{2,0}/(f)$, which shows
the claim. $\mathscr{B}_0^{\RR}(\gamma )= \ann_{\mathscr{B}^{\RR}_0}(f_2)$ is obviously. We have obtained
\[
\mathscr{C}_0^{\RR}=\frac{\mathds{B}}{\ann_{\mathds{B}}(c)}=\frac{\mathscr{E}_{2,0}}{(\gamma ,f)}.
\] 
To prove that Theorem \ref{them2} and Theorem \ref{gm} produce the same values one has to check,
that there is a positive real number $r$ so that $r c c_1=\det DX$ in $\mathds{B}$. We have
\[
c\cdot c_1=c\text{ trace }DX-c^2=c\det\frac{\partial (\gamma ,f)}{\partial (x,y)}
\]
in $\mathds{B}$. We verify the existence of such an $r$ in an example:

Set $f:=x^2+y^{k+1}$, $\gamma:=y$, $\delta_1:=-(k+1)$, $\delta_2:= y^l$. We have
$X_1=-(k+1)x^2$ and $X_2=-2xy+x^2y^l+y^{l+k+1}$. One computes
$\det DX= -2(k+1)(l+k+1)xy^{l+k}$ and $c\det\frac{\partial (\gamma ,f)}{\partial (x,y)}
=-2(k+1)xy^{l+k}$ in $\mathds{B}$. Here we have $r=l+k+1$.

\subsection*{Acknowledgements}
The author would like to thank W. Ebeling and X. G\'{o}mez-Mont for useful discussions about
indices of vector fields. He is also grateful to D. van Straten for examples which showed
that there was a mistake in the former proof of Theorem \ref{them1}. Finally he also wants to thank
the referees for useful and interesting comments on the paper.

\end{document}